\documentclass[11pt, reqno]{amsart}
 
 \usepackage{amsmath}
 \usepackage{amssymb}
\usepackage{bm}

\DeclareMathAccent{\mathring}{\mathalpha}{operators}{"17}

\newcommand{\mysection}[1]{\section{#1}
      \setcounter{equation}{0}}

\newtheorem{theorem}{Theorem}[section]
\newtheorem{lemma}[theorem]{Lemma}

\newtheorem{corollary}[theorem]{Corollary} 

\theoremstyle{definition}
\newtheorem{assumption}{Assumption}[section]

\theoremstyle{remark}
\newtheorem{remark}{Remark}[section]

\newcommand{\loc}{\text{\rm loc}}

 \makeatletter 
 \def\dashint{%
 \operatorname%
 {\,\,\text{\bf--}\kern-.98em\DOTSI\intop\ilimits@\!\!}}
\def\ninf{\qopname\relax\@empty{inf\phantom{p}\!\!\!}}

 \makeatother

\newcommand\bbeta{\text{\raise-.2ex\hbox{$\bm{\beta}$}}}

\newcommand\bR{\mathbb{R}}

\newcommand\bS{\mathbb{S}}

\def\cC{\mathcal{C}}

\newcommand\cF{\mathcal{F}}

\newcommand\cO{\mathcal{O}}

\def\sft{{\sf t}}

\begin{document}

\title[Diffusion
 processes with drift in a Morrey class]
{On   diffusion
 processes with drift in a Morrey class containing  $L_{d+2}$}

\author{N.V. Krylov}
 
\email{nkrylov@umn.edu}
\address{127 Vincent Hall, University of Minnesota,
 Minneapolis, MN, 55455}
 


\begin{abstract}
We present new conditions on the drift of the Morrey type
with mixed norms allowing us to obtain Aleksandrov type
estimates
of potentials of time inhomogeneous diffusion
processes in spaces with mixed norms and, for
instance, in $L_{d_{0}+1}$ with $d_{0}<d$. 
\end{abstract}

\maketitle

\mysection{Introduction and main results}

Let $\bR^{d}$ be a  Euclidean space of points
$x=(x^{1},...,x^{d})$, $d\geq 2$.  
Let $(\Omega,\cF,P)$ be a complete probability
space, let $\cF_{t}, t\geq0$, be an increasing family of
complete $\sigma$-fields $\cF_{t}\subset\cF$,  
and let $w_{t}$ be an $\bR^{d}$-valued Wiener process
relative to $\cF_{t}$. Fix $\delta\in (0,1)$
and denote by $\bS_{\delta}$
the set of $d\times d$ symmetric matrices
whose eigenvalues are between $\delta$ and
$\delta^{-1}$. 

\begin{assumption}
                                        \label{assumption 3.26.1}

On $\bR^{d+1}$
 we are given a smooth $\bS_{\delta}$-valued
function $\sigma(t,x)$ and a smooth
$\bR^{d}$-valued function $b(t,x)$ with compact support.
\end{assumption}

Under this assumption the solutions of the system
\begin{equation}
                                                 \label{11.29.2}
x _{s}=x  +\int_{0}^{s}\sigma (\sft_{r},x_{r})\,dw_{r}
+\int_{0}^{s}b (\sft_{r},x_{r}) \,dr,\quad \sft_{s}=t+s
\end{equation}
form a strong Markov process $(\sft_{s},x_{s})$. 
Our goal in this article is to
find conditions on the drift of the Morrey type
with mixed norms still allowing us to obtain Aleksandrov type
estimates
of potentials of $(\sft_{s},x_{s})$
 in spaces with mixed norms and, for
instance, in $L_{d_{0}+1}$ with $d_{0}<d$.

We take and fix in the whole article
two sets of numbers satisfying
\begin{equation}
                                              \label{3.21.70}
p_{0},q_{0}\in[1,\infty],\quad
  \frac{d_{0}}{p_{0}}+\frac{1}{q_{0}}=1,
\end{equation}
\begin{equation}
                                              \label{3.21.7}
p,q\in[1,\infty],\quad  \frac{d_{0}}{p}+\frac{1}{q}=1,
\end{equation}
where   $d_{0}=d_{0}(d,\delta)\in(d/2,d)$ is defined in
Section \ref{section 10.25.2}.
 
 Introduce
$$
B_{R} =\{x\in\bR^{d}:| x|<R\}, 	\quad B_{R}(x)=x+B_{R},
\quad C_{T,R}=[0,T)\times B_{R},
$$
$$
C_{T,R}(t,x)=(t,x)+C_{T,R},\quad C_{R}(t,x)=C_{R^{2},R}(t,x),
\quad C_{R} =C_{R}(0,0),
$$
and let
$\cC_{R}$ be the collection of cylinders $C_{R}(t,x)$,
$(t,x)\in\bR^{d+1}$,  $\cC=
 \{\cC_{R}:R>0\} $.

Introduce $\hat b $ as a constant
such that, for any $R\in(0,\infty)$ and $C\in\cC_{R}$,
\begin{equation}
                                                \label{3.26.20}
\|bI_{C}\|_{L_{p_{0},q_{0}} }\leq\hat b R^{d/p_{0}+2/q_{0}-1} ,
\end{equation}
where the norm $\|\cdot\|_{L_{p,q}}$ is introduced as follows.

For $p,q\in[1,\infty)$  
 we introduce the space $L_{p,q} $ as the space of Borel
functions on $\bR^{d+1}=\bR\times \bR^{d}$ such that
$$
\|f\|^{q}_{L_{p,q} }:=
\int_{\bR}\Big(\int_{\bR^{d}} |f(t,x)|^{p}\,dx\Big)^{q/p}\,dt<\infty 
$$
if $p\geq q$ or
$$
\|f\|^{p}_{L_{p,q} }:=\int_ {\bR^{d}} \Big(\int_{\bR} 
|f(t,x)|^{q}\,dt\Big)^{p/q}\,dx<\infty 
$$
if $ p\leq q$ with natural interpretation
of these definitions if $p=\infty$ or $q=\infty$.
We write $f\in L_{p,q}(Q)$ to mean that $fI_{Q}\in L_{p,q}$.
Observe that
$p $ is associated with   $x$ and
$q$ with   $t$ and the interior
integral is always elevated to the power $\leq 1$.  If $p=q$ we abbreviate $L_{p,p}$ to $L_{p}$, which is $L_{p}(\bR^{d+1})$. We use the same symbol $L_{p}$
for $L_{p}(\bR^{d})$ and hope that
its meaning will be clear from the context.
The necessity to define
$L_{p,q}$-norms differently for $p\geq q$
and $p\leq q$ is dictated by the form in which
we present our version
of the parabolic Aleksandrov estimate
  in Theorem \ref{theorem 3.27.2}.

Since $b$ is bounded and has compact support, such a
$\hat b<\infty$ satisfying \eqref{3.26.20} does exist.

We are now ready to present two main results of the article
in which
by   $\tau'_{R}$ we mean the first exit time
of $x_{s}$ from $B_{R}$. These theorems are proved in 
Sections \ref{section 4.2.1} and \ref{section 4.2.2}.

\begin{theorem}
                                          \label{theorem 3.27.2}
Under Assumption \ref{assumption 3.26.1}
there is a constant $\hat b<\infty$ depending only on
$d,\delta $, such that if \eqref{3.26.20}
holds for any $R\in(0,\infty)$ and $C\in\cC_{R}$, then
for   any $R\in(0,\infty)$, $ x
\in\bR^{d}$ 
and Borel $f\geq0$
\begin{equation}
                                              \label{3.30.3}
 I:= E_{0, x}\int_{0}^{\tau'_{R}}
f(t,x_{t})\,dt\leq
\hat N 
R^{2-d/p-2/q}\|f \|_{L_{p,q} },
\end{equation}
where $\hat N$ depends only on $d,\delta $.

\end{theorem}

It turns out that the global condition \eqref{3.26.20}
 can be replaced
with a local one at the expense of losing good control
on $\hat N$ in the estimate.

\begin{theorem}
                                          \label{theorem 3.27.20}
Let $R_{0}\in (0,\infty)$.
Under Assumption \ref{assumption 3.26.1}
there is a constant $\hat b<\infty$ depending only on
$d,\delta $, such that if 
\eqref{3.26.20} holds for any $R\in(0,R_{0}]$ and
$C\in \cC_{R}$, then
for   any $R\in(0,\infty)$, $ x
\in\bR^{d} $ 
and Borel $f\geq0$
\begin{equation}
                                              \label{3.21.80}
E_{0, x }\int_{0}^{\tau'_{R}}
f(t,x_{t})\,dt\leq
\hat N  \|f \|_{L_{p,q} },
\end{equation}
where $\hat N$ depends only on 
$d,\delta ,R$, and $R_{0}$.
\end{theorem}

\begin{remark}
                                             \label{remark 4.2.1}
By shifting the origin in $\bR^{d+1}$
one obtains estimates similar to \eqref{3.30.3}
and \eqref{3.21.80} for the process starting from any point  like
$$
E_{t,x}\int_{0}^{\tau'_{R}(y)}
f(\sft_{s},x_{s})\,ds\leq
\hat N  \|f \|_{L_{p,q} },
$$
where $\tau'_{R}(y)$ is the first exit time
of $x_{s}$ from $B_{R}(y)$.
 
\end{remark}

The above result have some implications
on relaxing  the  integrability requirement on $b$
in elliptic and parabolic Aleksandrov's
estimates.

In this part of the section Assumption \ref{assumption 3.26.1}
is replaced with the following.
\begin{assumption}
                                        \label{assumption 3.30.1}

On $\bR^{d+1}$
 we are given a Borel $\bS_{\delta}$-valued
function $\sigma(t,x)$ and a Borel
$\bR^{d}$-valued function $b(t,x)$
such that, for an $R_{0}\in(0,\infty)$,
estimate \eqref{3.26.20} holds with
$\hat b=\hat b(d,\delta)$ from Theorem \ref{theorem 3.27.2}
for any $R\in(0,R_{0}]$ and
$C\in \cC_{R}$.
\end{assumption}

Introduce
$$
 D_{i}=\frac{\partial}{\partial x^{i}},
\quad D_{ij}=D_{i}D_{j}\quad 
\partial_{t}=\frac{\partial}{\partial t},
$$
$$
a=(1/2)\sigma^{2},\quad
Lu=a^{ij}D_{ij}u+b^{i}D_{i}u.
$$
Recall that for a domain $Q\subset\bR^{d+1}$
one denotes by $\partial'Q$ its parabolic boundary
defined as the set of all points on $\partial Q$
which are endpoints of continuous curves of type $(t,x_{t})$,
$t\in[a,b]$, which
start in $Q$ and belong to $Q$ for all $t<b$.
Define $W^{1,2}_{p,q,loc}(Q)$ as the set of functions such that
$u,Du,D^{2}u,\partial_{t}u\in L_{p,q,loc}(Q)$.
Here is a quantitative maximum principle
that is a parabolic Aleksandrov estimate.

\begin{theorem}
                                         \label{theorem 3.30.1}
Under Assumption \ref{assumption 3.30.1}
 let $Q\subset\bR\times B_{R}$ be bounded and
 $u\in W^{1,2}_{p,q,\loc}(Q)\cap C(\bar Q)$.
  Take a function
$c\geq 0$ on $Q$. Then on $ Q$
\begin{equation}
                                               \label{10.14.10}
u \leq N 
\|I_{Q,u>0}( \partial_{t}u+Lu-cu)_{-}\|_{L_{p,q} }
+\sup_{\partial'Q}u_{+},
\end{equation}
where $N$ depends only on $d,\delta, R_{0}$,
and   $R$.

In particular (the maximum principle),
if $\partial_{t}u+Lu-cu \geq0$ in $Q$ and $u\leq0$
on $\partial'Q$, then $u\leq 0$ in $Q$.
\end{theorem}

The proof of this theorem coincides with that
of Theorem 5.1 of \cite{Kr_21_1} apart from
the point that after reducing the general case
to the one in which $u\in W^{1,2}_{p,q}(Q)$
and $b$ is bounded, here, we approximate $\sigma,b$, and
$u$ with smooth functions by using mollifiers.
Of course, we observe that the assumption
about $b$ is preserved under this operation.
After that, as in \cite{Kr_21_1}, It\^o's
formula and Theorem \ref{theorem 3.27.2}
allow us to finish the proof.

An adaptation of Theorem \ref{theorem 3.30.1}
to elliptic operators
yields an elliptic Aleksandrov estimate
and shows advantages of having mixed norm estimates
for parabolic operators.
Recall that we agreed to use the symbol $L_{p}$ to
 mean either $L_{p}(\bR^{d})$
or $L_{p}(\bR^{d+1})$. What is the meaning in each concrete
case will be quite clear from the context.
For instance, in the following theorem 
$L_{d_{0}}=L_{d_{0}}(\bR^{d})$. 

\begin{theorem}
                                         \label{theorem 4.1.1}
Under Assumption \ref{assumption 3.30.1}
let $\cO$ be a bounded domain in $\bR^{d}$,
$u\in W^{2}_{d_{0},\loc}(\cO)\cap C(\bar \cO)$.
Also assume that $a$ and $b$ are independent of $t$
and $p=p_{0}=d_{0}$ ($q=q_{0}=\infty$).
Take a function $c\geq0$. Then on $\cO$
\begin{equation}
                                               \label{4.1.5}
u \leq N 
\|I_{\cO,u>0}(  Lu-cu)_{-}\|_{L_{d_{0}} }
+\sup_{\partial\cO}u_{+},
\end{equation}
where $N$ depends only on $d,\delta, R_{0}$, and
the diameter of
$ \cO$.

\end{theorem}

Proof. Take $\varepsilon>0$,
let $v(t,x)=e^{-\varepsilon t}u(x)$ and apply \eqref{10.14.10}
to $v$ and $(0,T)\times \cO$ in place of $u$ and $Q$,
respectively.
Then we get that on $\cO$
$$
u\leq N \|I_{\cO,u>0}(-\varepsilon u+  Lu-cu)_{-}\|_{L_{d_{0}} }
+\sup_{\partial\cO}u_{+}+e^{-T}\sup_{\cO}u_{+}.
$$
Letting $T\to\infty$ eliminates the last term,
After that letting $\varepsilon\downarrow 0$ 
 and observing that
$(-\varepsilon u+  Lu-cu)_{-}\leq \varepsilon
\sup|u|+( Lu-cu)_{-}$
yields 
\eqref{4.1.5}. The theorem is proved.

In case when $d_{0}$ in Theorem \ref{theorem 4.1.1}
is replaced by $d$ and $b\in L_{d}$, the result
  belongs to
A.D. Aleksandrov (1960), see Theorem 8 in
\cite{Al_60}.
The proofs are given in 1963 in \cite{Al_63}.
Our Theorem \ref{theorem 4.1.1} extends Aleksandrov's
result in reducing the power of summability of both:
the drift and the free terms. However, we treat only the
uniformly nondegenerate case.

 There was a
considerable interest in reducing $L_{d}$-norm
of the free term to $L_{d_{0}}$-norm with
$d_{0}<d$. This was achieved by Cabr\'e
\cite{Ca_95} for bounded $b$
and by Fok in \cite{Fo_98} for
  $b\in L_{d+\varepsilon}$.
In \cite{Kr_19_1} the  author
allowed $b\in L_{d}$ and the free term  in
$L_{d_{0}}$   with
$d_{0}<d$. This made it possible
to develop  in  \cite{Kr_21b}   a $W^{2}_{p}$-solvability
theory for linear equations with $b\in L_{d}$ and $p<d$. Applied
to fully nonlinear equations we can now
treat $W^{2}_{d_{0}}$-solvability with
``the coefficients'' of the first order terms
in $L_{d}$ (see \cite{Kr_20}).
 In Theorem 1.1 of \cite{DK_21} our
Theorem \ref{theorem 4.1.1} is proved, loosely speaking,
when $p=p_{0}=d_{0}$, but $R_{0}$ is comparable with the
diameter of $\cO$. By the way, observe that, if
$|b(x)|=c/|x|$ with small enough $c$, then $b$
 satisfies the assumption
of Theorem \ref{theorem 4.1.1}
(and of Theorem 1.1 of \cite{DK_21}) and is not of class
$L_{d,\loc}$.

The results like Theorem \ref{theorem 3.30.1} are
indispensable in the theory of
controlled diffusion processes (see, for instance,
\cite{Kr_77}). First such  result with bounded $b$ and
$d$ in place of $d_{0}$ and $L_{d+1}$ in place of 
$L_{p,q}$ was published in \cite{Kr_74}. 
It  was extended by A.I. Nazarov and N.N. Ural'tseva
\cite{NU_85} to allow $b\in L_{d+1}$. The author in
\cite{Kr_86} developed a general approach to such estimates
and slightly improved the result of \cite{NU_85}.
By using this approach A.I. Nazarov in
\cite{Na_15}
for the first time proved estimates
in $L_{p,q}$-spaces with $b\in L_{p_{0},q_{0}}$
when $d_{0}$ in \eqref{3.21.70} and \eqref{3.21.7}
is replaced by $d$ but no further restrictions like
\eqref{3.26.20} are imposed.

 In \cite{CKS_00} by extending some earlier
results by Wang the authors prove Theorem \ref{theorem 3.30.1}
for $L_{p}$-viscosity solutions with
$p=q<d+1$   when $b$ is   bounded.
 In our situation
we have some freedom in choosing $p,q$ and $b\in L_{p_{0},q_{0}}$,
but we only treat true solutions. Our
Theorem \ref{theorem 3.30.1} covers
Theorem 2.4 of \cite{CKS_00} on the account of
having mixed norms
and $b\in L_{p_{0},q_{0}}$.

The author in \cite{Kr_21_1} gave a version of the result
of \cite{Na_15} by allowing $f\in L_{p,q}$
with $p,q$ satisfying \eqref{3.21.7} and
$b\in L_{p_{0},q_{0}}$ but $p_{0},q_{0}$
are supposed to satisfy \eqref{3.21.70}
with $d$ in place of $d_{0}$ {\em and\/}
\eqref{3.26.20} is imposed. 
Relaxing the assumptions on $f$ further
in \cite{Kr_21_2} allowed the author
to investigate fine properties of the corresponding
diffusion processes and lead in \cite{Kr_21_3}
to proving the solvability in Sobolev space  $W^{1,2}_{p}$,
where $p<d+1$,
of the equation
\begin{equation}
                                      \label{12.23.1}
\partial_{t}u+\Delta u+b^{i}D_{i}u=f
\end{equation}
when $f\in L_{q}$ with $q$ large enough, $b\in L_{d+1}$,
and \eqref{3.26.20} is satisfied with $p_{0}=q_{0}=d+1$.
In this case it turns out that $b^{i}D_{i}u\in L_{p}$.

In this article we relax the restrictions on $b$
from \cite{Kr_21_1} allowing $p_{0},q_{0}$  to be
as in \eqref{3.21.70}. Observe that the function $b$ with
$$
|b(t,x)|=I_{C_{1}}(t,x)\frac{1}{|x|+\sqrt t}\Big(\frac{|x|}{
\sqrt t}\Big)^{2/(d+1)}
$$
satisfies \eqref{3.26.20} with finite $\hat b$ and $p_{0}=q_{0}=
d_{0}+1$ but does not belong to $L_{d+1}$.

Condition \eqref{3.26.20},
however, does not allow us to derive from the results
presented here even the estimates
from \cite{NU_85} and \cite{Kr_86}. The point is that,
if we have $b\in L_{r}$ and ask ourselves what $r$
should be in order to have $b\in L_{p_{0},q_{0}}$
{\em satisfying\/} \eqref{3.26.20}, then the answer
is: $r\geq d+2$
 (thus explaining the title
of the paper). So $b\in L_{d+1}$ are not good enough.
This is somewhat discouraging but, most likely
$b\in L_{d+1}$ is not good enough to have
solvability in any $W^{1,2}_{r}$ of \eqref{12.23.1}.
In any case, by imposing \eqref{3.26.20}, say with $p_{0}=q_{0}=d_{0}+1$,
we allow stronger local singularities 
of $b$, as compared to $b\in L_{d+1}$, spread sufficiently far apart.

This article have some similarity to \cite{Kr_21_1},
in particular, we borrow from there quite a few results.
The main difference is that in \cite{Kr_21_1}
we start by considering general It\^o processes,
which led to the requirement $q_{0}<\infty$,
and here we start by considering Markov diffusion
processes with regular drift and diffusion coefficients.
Therefore we narrow the class of processes under
investigation but gain obtaining better estimates
by exploiting the fact that we have solutions
of equations coefficients of which we can change
(see Section \ref{section 4.2.1}).
Our results could    possibly  be relevant
in investigations described in \cite{BFGM_19},
\cite{RZ_21_1}, \cite{RZ_21_2} and the references therein,  where the authors present very strong results on fine
properties of diffusion processes under
some regularity assumptions on $\sigma$
and $b\in L_{p,q}$ with $d/p+2/q\leq 1$. We add
some information about Green's functions
of such processes.

 The rest of the article is organized as follows.
Section \ref{section 10.25.2} contains some auxiliary
results part of which is borrowed from \cite{Kr_21_1}.
The other part of the section is devoted
to proving better summability of Green's functions
and mixed norm estimates.
In Section \ref{section 4.2.1}
we prove Theorem \ref{theorem 3.27.2}.
Theorem \ref{theorem 3.27.20} is proved in Section
\ref{section 4.2.2}. Section \ref{section 3.21.1}
is an Appendix where we prove a version of Gehring's
lemma used in Section \ref{section 10.25.2}.

We finish the section by some notation.  
If $\Gamma$ is a measurable subset of $\bR^{d+1}$ we denote by
$|\Gamma|$ its Lebesgue measure. The same notation
is used for measurable  subsets of $\bR^{d}$ with $d$-dimensional
Lebesgue measure. We hope that it will be clear
from the context which Lebesgue measure is used.
If $\Gamma$ is a measurable subset of $\bR^{d+1}$ and
$f$ is a function on $\Gamma$ we denote 
$$
\dashint_{\Gamma}f\,dxdt=\frac{1}{|\Gamma|}
\int_{\Gamma}f\,dxdt.
$$
In case $f$ is a function on a measurable subset $\Gamma$
of $\bR^{d}$ we set
$$
\dashint_{\Gamma}f\,dx =\frac{1}{|\Gamma|}
\int_{\Gamma}f\,dx .
$$

\mysection{Auxiliary results}

                                        \label{section 10.25.2}

In this section we do not suppose that any version of
\eqref{3.26.20} holds.

For $C\in \cC_{R}$
denote by $\tau_{C}$ the first exit time of
$(\sft_{s},x_{s})$ from $C$, but if $C=C_{R}$,
we use the notations $\tau_{R}$ instead of $\tau_{C_{R}}$
with the hope that no confusion will be created. 

Here is a combination of Lemmas 4.1 and 4.2 of \cite{Kr_20_2}.

\begin{theorem}
                                        \label{theorem 9.27.1}

For any $(t,x)\in\bR^{d+1}$, Borel $f(s,y),g(y)\geq0$ and stopping time $\gamma$
\begin{equation}
                                                   \label{5.6.4}
  E_{t,x}\int_{0}^{\gamma}   
f( \sft_{t},x_{t})\,dt\leq N(d, \delta  )
\big(A+ B^{2}
\big)^{d/(2d+2 )}
\|f\|_{L_{d+1 }},
\end{equation}
\begin{equation}
                                                   \label{3.7.1}
  E_{t,x}\int_{0}^{\gamma}   
g(  x_{t})\,dt\leq N(d, \delta  )
\big(A+ B^{2}
\big)^{1/2}
\|g\|_{L_{d  } },
\end{equation}
where $A=E_{t,x}\gamma$ and
$$
B=E_{t,x}\int_{0}^{\gamma}|b(\sft_{s},x_{s})|\,ds.
$$

\end{theorem}
 
Observe that if $\gamma=\tau_{C}$ with $C\in\cC_{R}$
in Theorem \ref{theorem 9.27.1}, 
then obviously $\gamma
\leq R^{2}$ and $A\leq R^{2}$. Since $b$ is bounded,
  $B\le KR^{2}$, where $K$ is independent of $C$. This shows that
$$
\bar b_{R}:=\sup_{\substack{\rho\leq R \\ \,C\in\cC_{\rho}}}
\frac{1}{\rho}
\sup_{(t,x)\in\bR^{d+1}}E_{t,x}\int_{0}^{\tau_{C}}
|b(\sft_{s},x_{s})|\,ds
$$
is finite for any $R\in(0,\infty)$.
 
\begin{remark}
                                       \label{remark 3.27.1}
Usual way to deal with additive functionals shows that
for any $n=1,2,...$, $R\in(0,\infty)$, $(t,x)\in\bR^{d+1}$,
and $C\in\cC_{R}$
$$
E_{t,x}\Big(\int_{0}^{\tau_{C}}
|b(\sft_{s},x_{s})|\,ds\Big)^{n}\leq n!\,\bar b_{R}^{n}R^{n}.
$$

\end{remark}

 Since $B\leq \bar b_{R} R$ if $\gamma=\tau_{C}$
and $C\in\cC_{R}$, we obtain the following.

\begin{lemma}
                                      \label{lemma 8.16.1}
 
For any Borel $f \geq0$, $R>0$,
$C\in\cC_{R}$, and $(t,x)\in\bR^{d+1}$
we have
\begin{equation}
                                          \label{9.29.2}
 E_{t,x}\int_{0}^{\tau_{C} }  
f( \sft_{s},x_{s})\,ds\leq 
 N(d,  \delta)(1+\bar b_{R}) ^{d/  (d+1)  }
 R 
 ^{d/ (d+1)  }
\|f\|_{L_{d+1 }}.
\end{equation}
 
\end{lemma}

By plugging in $f=|b|I_{C}$ and using the fact that $b$
is bounded and has compact support, we get
$$
E_{t,x}\int_{0}^{\tau_{C} }  
|b( \sft_{s},x_{s})|\,ds
$$
$$
\leq 
K(1+\bar b_{R}) ^{d/  (d+1)  }
 R 
 ^{d/ (d+1)  }
\|f\|_{L_{d+1 }}\leq K(1+\bar b_{R}) ^{d/  (d+1)  }R,
$$
where the constants $K$ are independent 
of $R>0$,
$C\in\cC_{R}$, and $(t,x)\in\bR^{d+1}$. It 
follows by definition that
$$
\bar b_{R}\leq K(1+\bar b_{R}) ^{d/  (d+1)  },
$$
which shows that $\bar b_{\infty}:=\lim_{R\to\infty}
\bar b_{R}$ is finite.  

 For $T,R\in(0,\infty)$ introduce
$$
\tau_{T,R} =\inf\{s\geq 0:(\sft_{s}, x_{s})\not\in C_{T,R}\} 
\quad (\tau_{R}=\tau_{R^{2},R}) .
$$
  
  Observe that owing to \eqref{9.29.2} with $f=I_{C}$
\begin{equation}
                                          \label{3.29.1}
E_{0,x}\tau_{R}\leq N(d,\delta)(1+\bar b_{\infty})^{d/(d+1)}R^{2}.
\end{equation}

Estimate \eqref{3.29.1} says that   $\tau_{R}$ is of order
not more than $R^{2}$.
A very important fact which is implied by Corollary
\ref{corollary 7.29.1} is that 
$\tau_{R}$ is of order
not less than $R^{2}$.
To show this we need
an additional assumption appearing after the following result,
in which
\begin{equation}
                                                  \label{1.6.1}
\tau'_{R} =\inf\{t\geq0: x_{t}\not\in B_{R} \},\quad
\gamma_{R} =\inf\{t\geq0:x_{t} \in \bar B_{R} \}.
\end{equation}

\begin{theorem}
                                      \label{theorem 8.2.1}
There are
   constants $\bar \xi=\bar \xi(d,\delta)\in (0,1) $ and
 $\bar N=\bar N(d,  \delta)$ {\em
continuously\/} depending on $\delta$
such that if, for an $  R\in(0,\infty)$, we have
\begin{equation}
                                     \label{12.18.3}
\bar N \bar b_{  R}\leq 1,
\end{equation}
then  for any $x$
\begin{equation}
                                          \label{8.2.2} 
  P_{0,x}( \tau_{R}  =   R^{2} )\leq 1-\bar\xi,
\quad P_{ 0,0}( \tau_{R}  =   R^{2} )\geq 
 \bar\xi .   
\end{equation}
Moreover for $n=1,2,...$  
\begin{equation}
                                          \label{1.3.1} 
P_{0,x}(\tau'_{R}(x)\geq nR^{2})
=  P_{0,x}( \tau_{nR^{2},R}(x)  =  n R^{2} )\leq (1-\bar\xi)^{n},   
\end{equation}
so that 
\begin{equation}
                                          \label{3.7.2}
E_{0,x}\tau'_{R} \leq N(d,\delta)R^{2},
\end{equation}
and
\begin{equation}
                                          \label{1.3.3}  
I:=E_{0,x}\int_{0}^{\tau'_{R} }|b(t,x_{t})|\,dt
\leq N(d ,\delta) \bar b_{R}R.
\end{equation}

Furthermore,  the probability
starting from a point in  $\bar B_{9R/16}$
to reach the ball $\bar B_{R/16}$ before exiting
from $B_{R}$ is bigger than $\bar\xi$:
for any $x $ with $|x|\leq 9R/16$ 
\begin{equation}
                                          \label{1.2.1} 
P_{0,x}(\tau'_{R} >\gamma_{R/16} )\geq\bar\xi.   
\end{equation}

\end{theorem}

Proof. This theorem is similar to Theorem 2.3 of \cite{Kr_21_1}.
The most significant difference  is that 
in Theorem 2.3 of \cite{Kr_21_1}
to estimate quantities like
\begin{equation}
                                          \label{3.1.1} 
E_{0,0}\int_{0}^{\tau_{R} }  
|b(s,x_{s})|\,ds
\end{equation}
an estimate similar to \eqref{9.29.2} is used and this lead
to the assumption that $q_{0}<\infty$. In our situation
quantity \eqref{3.1.1}  is less than $\bar b_{R}R$ by definition.
Taking this into account to prove the theorem
one can just repeat the proof of Theorem 2.3 of \cite{Kr_21_1}.

We show an example of how to do that
proving  \eqref{1.3.3}. By using the strong Markov property
we obtain
$$
I=\sum_{n=1}^{\infty}E_{0,x}I_{\tau_{(n-1)R^{2},R} 
>\tau'_{R} }E\Big(
\int_{\tau_{(n-1)R^{2},R} }^{\tau_{nR^{2},R} }
|b(t,x_{t})|\,dt\mid \cF_{\tau_{(n-1)R^{2},R}}\Big)
$$
$$
\leq J\sum_{n=1}^{\infty}P_{0,x}(\tau_{(n-1)R^{2},R}  
=(n-1)R^{2})\leq J\sum_{n=1}^{\infty}
(1-\bar\xi)^{n-1},
$$
where $J=R\bar b_{R}$.
This yields \eqref{1.3.3}. The theorem is proved.

In light of \eqref{3.7.2} and \eqref{1.3.3} 
estimate \eqref{3.7.1} implies the following.

\begin{corollary}
                                     \label{corollary 3.7.1}
If \eqref{12.18.3} holds, then    
for any Borel $g\geq0$ 
  and $x\in\bR^{d}$
we have
\begin{equation}
                                          \label{9.29.20}
 E_{0,x}\int_{0}^{\tau'_{R} }  
g( x_{s})\,ds\leq 
 N(d,  \delta)(1+\bar b_{R})  
 R 
\|g\|_{L_{d  }}.
\end{equation}
\end{corollary}

\begin{assumption}
                                   \label{assumption 12.18.2}
It holds  that
\begin{equation}
                                     \label{1.23.2}
\bar N(d,\delta) \bar b_{ \infty}< 1.
\end{equation}
\end{assumption}

This assumption (as well as Assumption  
\ref{assumption 3.26.1})
is supposed to hold throughout the section.

We need a few more results given here without proofs
because their proofs are obtained just by
repeating the corresponding proofs from
\cite{Kr_21_1}. 

Here are particular cases of Theorem 2.6 and
Corollary 2.7 of \cite{Kr_21_1}.

\begin{theorem}
                                     \label{theorem 8.20.1}
For any $\lambda,R>0$,   we have
\begin{equation}
                                          \label{8.20.1}
E_{0,0}e^{-\lambda  \tau_{R} }\leq
e^{\bar\xi/2}e^{- \sqrt{\lambda}  R
\bar\xi/2}.
\end{equation}
In particular, for  
  any   $t,R>0$   we have
\begin{equation}
                                             \label{10.2.2}
P_{0,0}(  \tau_{R}  \leq  t )\leq 
 e^{\bar\xi/2}\exp\Big(-\frac{{\bar \xi}^{2}R^{2}}{16 t}\Big).
\end{equation}

\end{theorem}

\begin{corollary}
                                       \label{corollary 7.29.1}
 
There is a  constant   
  $N =N (  \bar\xi)$
such that for any $R\in(0,\infty)$ 
\begin{equation}
                                                   \label{8.21.1}
N E _{0,0}\tau_{R} \geq R^{2} .
\end{equation}

\end{corollary}

A few more general results are related to going through
a long ``sausage".
\begin{theorem}
                                        \label{theorem 1.24.1}
Let $R\in(0,\infty)$, $x,y\in\bR^{d}$ and $16|x-y|\geq 3R$.  
For $r>0$ denote by $S_{r}(x,y)$   the open convex hull
of $B_{r}(x)\cup B_{r}(y)$. Then there exist
$T_{0},T_{1}$, depending only on $\bar\xi$,
such that $0<T_{0}<T_{1}<\infty$ and the $P_{0,x}$-probability $\pi$
that $x_{t}$ will reach $\bar B_{R/16}(y)$ before exiting
from $S_{R}(x,y)$ and this will happen
on the time interval $[nT_{0}R^{2},nT_{1}R^{2}]$
is greater than $\pi_{0}^{n}$, where
$$
n= \Big\lfloor \frac{16|x-y|+R}{4R}\Big\rfloor 
$$  
and $\pi_{0}=\bar\xi/3$.

\end{theorem}  

This is a particular case of Theorem 2.9 of \cite{Kr_21_1}.

\begin{remark}
                                      \label{remark 1.25.1}
Observe that, for any fixed $x,y$, the interval
$[nT_{0}R^{2},nT_{1}R^{2}]$ is as close to zero
as we wish if we choose $R$ small enough.
Then, of course, the corresponding probability will be
quite small but $>0$.

\end{remark} 

Next follow analogs of Corollaries 2.10 and 2.11 of
\cite{Kr_21_1}
\begin{corollary}
                                     \label{corollary 1.25.1}
Let   $\kappa\in[0,1)$, and $|x|\leq\kappa R$.
Then for any $T>0$
\begin{equation}
                                        \label{1.25.2}
NP_{0,x}(\tau'_{R}> T)\geq e^{-\nu T/[(1-\kappa)R]^{2}},
\end{equation}
where $N$ and $\nu>0$ depend only on $\bar\xi$.
\end{corollary}

\begin{corollary}
                                  \label{corollary 2.3.1}
Let   $R\in(0,\infty) $. Then there
exists a constant $N$, depending only on
$\bar\xi $, such that, for any 
$T>0$,
$$
P_{0,0}(\tau'_{R}>T)\leq Ne^{-T/(NR^{2})}.
$$
\end{corollary}

It is well known that, in light of
Assumption \ref{assumption 3.26.1},
the process $(\sft_{s},x_{s})$ has a Green's function,
which means that we can 
introduce a function   $G (t,x,s,y)\geq0$ so that   for all nonnegative Borel
$f$ and $(t,x)\in\bR^{d+1}$
$$
E_{t,x}\int_{0}^{\infty}f(\sft_{r},x_{r}
)\,dr=\int_{t}^{\infty}\int_{\bR^{d } }G (t,x,s,y)
f(s,y)\,dyds,
$$
$$
 G(s,y)=G (0,0,s,y).
$$

\begin{theorem}
                                           \label{theorem 9.3.2}
There exist  
$d_{0}\in(1,d)$, depending only
on $\delta,d $,   and a constant    $N=N(\delta,d) $
such that for any     $R\in(0,\infty)$, $C\in\cC_{R}$,
 and $p\geq d_{0}+1 $, we have
\begin{equation}
                                          \label{10.14.1}
\Big(\dashint_{C }G^{p/(p-1)}(s,y)
\,dyds\Big)^{(p-1)/p}
\leq N  R^{-d} .
\end{equation}
 
\end{theorem}

Proof. We basically follow the idea in \cite{FS_84}.
Introduce $\cC_{+}$ as the set of cylinders $C_{R}(t,x)$, $R>0$,  
$t\geq0$, $x\in\bR^{d}$. For $C=C_{R}(t,x)\in\cC_{+}$ let
 $2C=C_{2R}(t,x)$.    If $C=C_{R}(t,x)$ set $R_{C}=R$.

Take $C\in\cC_{+}$   and
define 
recursively 
$$
\gamma^{1}=\inf\{s\geq0 :(\sft_{s},x_{s})\in \bar C\},\quad
\tau^{1}=\inf\{s\geq\gamma^{1} :(\sft_{s},x_{s})\not\in 2C\},
$$
$$
\gamma^{n+1}=\inf\{s\geq\tau^{n} :(\sft_{s},x_{s})\in \bar C\},\quad
\tau^{n+1}=\inf\{s\geq\gamma^{n+1} :(\sft_{s},x_{s})\not\in 2C\}.
$$
Then for any nonnegative Borel $f$ vanishing outside $C$
with $\|f\|_{L_{d+1}(C)}=1$
we have
$$ 
\int_{C}f (s,y)  G(s,y)\,dyds
=E_{0,0}\int_{0}^{\infty } f (\sft_{s},x_{s})\,ds 
$$
$$
=\sum_{n=1}^{\infty}E_{0,0} 
E \Big(\int_{\gamma^{n}}^{\tau^{n}} f (\sft_{s},x_{s})\,ds
\mid \cF_{\gamma^{n}}\Big).
$$
By the strong Markov property and \eqref{9.29.2} (a.s.)
($\bar b_{\infty}$ is estimated by using \eqref{1.23.2})
$$
E\Big(\int_{\gamma^{n}}^{\tau^{n}} f (\sft_{s},x_{s})\,ds
\mid \cF_{\gamma^{n}}\Big)\leq
N(d,  \delta) 
 R ^{d/ (d+1)  }_{C}I_{\gamma^{n}<\tau^{n}},
$$
Next, on the set $\{\gamma^{n}<\tau^{n}\}$ we have
$(\sft_{\gamma_{n}},x_{\gamma_{n}})\in \bar C$
and by Corollary \ref{corollary 7.29.1} on average
it will take at least $R_{C}^{2}/N$ time for the process
to exit from $2C$, that is
$$
E(\tau^{n}-\gamma^{n}\mid \cF_{\gamma^{n}})
\geq R_{C}^{2}/N
$$
on $\{\gamma^{n}<\tau^{n}\}$. We conclude that
$$ 
\int_{C}f (s,y)  G ( s,y)\,dyds
\leq N R_{C} ^{-(d+2)/ (d+1)}
\sum_{n=1}^{\infty}E_{0,0}(\tau^{n}-\gamma^{n})
$$
$$
\leq N R _{C}^{-(d+2)/ (d+1)}
E_{0,0}\int_{0}^{\tau_{C_{0}}}I_{2C}(\sft_{s},x_{s})\,ds
$$
$$
=
N R _{C}^{-(d+2)/ (d+1)}
\int_{2C}G ( s,y)\,dyds.
$$

The arbitrariness of $f$ 
yields 
$$
\Big(\dashint_{C}G^{(d+1)/d} \,dyds\Big)^{d/(d+1)}
\leq N  
\dashint_{2C }G \,dyds.
$$
Now, by Theorem 
\ref{theorem 3.9.1} there exists $d_{0}=d_{0}(d,\delta)<d$
such that 
$$
\Big(\dashint_{C}G^{(d_{0}+1)/d_{0}} 
\,dyds\Big)^{d_{0}/(d_{0}+1)}
\leq N 
\dashint_{2C }G \,dyds
$$
for any $C\in \cC_{+}$. H\"older's inequality shows that this
estimate also holds for $p\geq d_{0}$ in place of $d_{0}$.
After that it only remains to recall that
 the integral of $G$ over $[0,R^{2}]\times \bR^{d}$
is $R^{2}$ and hence
$$
\dashint_{2C }G \,dxdt\leq N R^{-d}_{C}.
$$
 The theorem is proved.

\begin{theorem}
                      \label{theorem 9.1.1}
For   any $R\in(0,\infty)$, $ x
\in\bR^{d} $ 
and Borel $f\geq0$
\begin{equation} 
                            \label{9.1.2}
E_{0,x}\int_{0}^{\tau'_{R}}
f(t,x_{t})\,dt\leq
\hat N 
R^{(2d_{0}-d )/(d_{0}+1)}\|f \|_{L_{d_{0}+1} },
\end{equation}
where $\hat N$ depends only on $d,\delta $.

\end{theorem}

Proof. Observe that in light of Theorem
\ref{theorem 9.3.2} and the   Markov property
$$
E_{0,x}\int_{0}^{\tau'_{R}}
f(s,x_{s})\,dt=\sum_{n=0}^{\infty}
E_{0,x}I_{\tau'_{R}\geq nR^{2}}
E\Big(\int_{nR^{2}}^{((n+1)R^{2})\wedge\tau'_{R}}f(\sft_{s},x_{s})\,ds\mid 
\cF_{nR^{2}}\Big)
$$
$$
\leq N\|f\|_{L_{d_{0}+1}}R^{(d+2)d_{0}/(d_{0}+1)-d}\sum_{n=0}^{\infty}
P_{0,x}(\tau'_{R}\geq nR^{2}).
$$
It only remains to use Corollary 
\ref{corollary 2.3.1}. The theorem is proved.

Introduce the Green's function $G_{R}(x,y)$ of $x_{\cdot}$
 by means of the formula
$$
E_{0,x}\int_{0}^{\tau'_{R}}f(x_{t})\,dt=\int_{B_{R}}
f(y)G_{R}(x,y)\,dy, .
$$
In this ``elliptic'' setting one can use Corollary
 \ref{corollary 3.7.1} in place of \eqref{9.29.2}
and then by mimicking the proof of Theorem
\ref{theorem 9.3.2} one sees that there 
is a constant $N=N(d,\delta)$ such that
for any $R\in(0,\infty)$ and any ball $B$
such that $4B\subset B_{R}$ it holds that
$$
\Big(\dashint_{B}G_{R}^{d/(d-1)}( x ,y)\,dy\Big)^{(d-1)/d}
\leq N  \dashint_{2B}
G_{R}(x,y)\,dy.
$$
Then by following the arguments in \cite{GM_79}
(or following our arguments in Section \ref{section 3.21.1})
one can see that
there exists $d_{0}=d_{0}(d,\delta)<d$ such that
\begin{equation}
                                               \label{3.21.6}
\Big(\dashint_{B}G_{R}^{p/(p-1)}(x,y)\,dy\Big)^{(p-1)/p}
\leq N  \dashint_{2B}
G_{R}(x,y)\,dy
\end{equation}
for any $p\geq d_{0}$ and $B$
such that $4B\subset B_{R}$. Of course, we can take
the $d_{0}$'s to be the same.

Next \eqref{3.21.6} implies that
$$
\Big(\dashint_{B_{R/4}}G_{R}^{p/(p-1)}(x,y)\,dy\Big)^{(p-1)/p}
\leq N  \dashint_{B_{R}}
G_{R}(x,y)\,dy
\leq N R^{2-d},
$$
where the second inequality follows from \eqref{3.7.2}.
In the inequality between the extreme terms one
can replace $G_{R}$ by a smaller quantity $G_{R/4}$.
Then by using the arbitrariness of $R$
 we arrive at the following.

\begin{theorem}
                                      \label{theorem 3.21.1}
For $p\geq d_{0}$ and any $R\in(0,\infty)$
and $ x\in \bR^{d} $
$$
\Big(\dashint_{B_{R }}G_{R}^{p/(p-1)}(x,y)\,dy\Big)^{(p-1)/p}
 \leq N (d,\delta) R^{2-d}.
$$
\end{theorem}

\begin{remark}
If in Theorem \ref{theorem 3.21.1} we take $x_{t}$
to be just a Wiener process, then we will see that 
$d_{0}>d/2$.
\end{remark}

In the following theorem we use the interpolation technique
suggested by A.I. Nazarov in \cite{Na_15}.

\begin{theorem}
                                      \label{theorem 3.21.2}
There exists a constant $\hat N=\hat N(d,\delta )$
such  that for any $R\in(0,\infty)$,
$ x\in \bR^{d} $ and Borel $f\geq0$  
estimate \eqref{3.30.3} holds.

\end{theorem}

Proof.   If $p=d_{0}+1$, then $q=d_{0}+1$  and 
\eqref{3.30.3} follows from Theorem \ref{theorem 9.1.1}  
since   $(2d_{0}-d )/(d_{0}+1)
=2-d/(d_{0}+1)-2/(d_{0}+1)$.  
 
If $p=d_{0}$ and $q=\infty$ estimate \eqref{3.30.3}
follows from Theorem \ref{theorem 3.21.1} since
$$
I\leq E_{0,x}\int_{0}^{\tau'_{R}}\sup_{s\geq0}f(s,x_{t})
\,dt=\int_{B_{R}}G_{R}(x,y)\sup_{s\geq0}f(s,y )\,dy
$$
$$
\leq \|f\|_{L_{d_{0},\infty} }  
\|G_{R}\|_{L_{d_{0}/(d_{0}-1)}
(B_{R})} \leq N  R^{2-d/d_{0}}
\|f\|_{L_{d_{0},\infty} }.
$$
If $p=\infty$ and $q=1$,   estimate \eqref{3.30.3}
 holds because
$$
I\leq E_{0, x}\int_{0}^{ \infty}\sup_{\bR^{d}}f(t,y)
\,dt=\|f\|_{L_{p,q} }.
$$

In case $\infty>p>d_{0}+1$ we have $p>q$ and
set $\beta=p/(d_{0}+1)$ and $\alpha=\beta/(\beta-1)$.
Take a nonnegative $g(t)$ such that $ \big(f(t,x)g(t)
\big)/g(t)=f(t,x)$ ($0/0=0$) and use H\"older's inequality to conclude
that $I\leq I_{1}I_{2}$, where
$$
I_{1}=\Big(\int_{0}^{ \infty} g^{-\alpha}(t)
\,dt\Big)^{1/\alpha},
$$
$$
I_{2}=\Big(E_{0,x}\int_{0}^{ \tau'_{R}} 
g^{\beta}(t)f^{\beta}(t,x_{t})
\,dt\Big)^{1/\beta}
$$
$$
\leq N 
R^{(2d_{0}-d)/p}
 \Big(\int_{0}^{ \infty}g^{(d_{0}+1)\beta}(t)
\Big(\int_{\bR^d}f^{(d_{0}+1)\beta}(t,y)\,dy\Big)\,dt
\Big)^{1/(d_{0}\beta+\beta)}.
$$
For $g$ found from
$$
  g^{-\alpha}(t)=
g^{(d_{0}+1)\beta}(t)
\int_{\bR^d}f^{(d_{0}+1)\beta}(t,y)\,dy
$$
we get \eqref{3.30.3} and this takes care of the case
that $\infty>p>d_{0}+1$.

If $\infty>q>d_{0}+1$ we have $p<q$ and
set $\beta=q/(q-d_{0}-1)$ and $\alpha=\beta/(\beta-1)$.
Take a nonnegative $g(y)$ such that $ \big(f(t,y)g(y)
\big)/g(y)=f(t,y)$ ($0/0=0$) and use H\"older's 
inequality to conclude
that $I\leq I_{1}I_{2}$, where
$$
I_{1}=\Big(E_{0,x }\int_{0}^{\tau_{R}}
g^{-\beta}(x_{t})\,dt
\Big)^{1/\beta}
\leq N R^{(2-d/d_{0})/\beta}
\Big(\int_{\bR^d}g^{-d_{0}\beta}(y)\,dy
\Big)^{1/(d_{0}\beta)},
$$
$$
I_{2}=\Big(E_{0,x }\int_{0}^{ \tau'_{R}}
g^{\alpha}(x_{t})
f^{\alpha}(t,x_{t})\,dt
\Big)^{1/\alpha}
$$
$$
\leq N 
R^{(2d_{0}-d)/(\alpha d_{0}+\alpha)}
\Big(\int_{\bR^{d}}g^{(d_{0}+1)\alpha}(y)
\Big(\int_{0}^{\infty}
f^{(d_{0}+1)\alpha}(t,y)\,dt\Big)dy
\Big)^{1/(\alpha d_{0}+\alpha )}.
$$
For $g$ found from
$$
g^{-d_{0}\beta}(y)=g^{(d_{0}+1)\alpha}(y)
 \int_{0}^{\infty}
f^{(d_{0}+1)\alpha}(t,y)\,dt 
$$
we get \eqref{3.30.3} after simple manipulations
and this proves the theorem.

Here is a key to proving Theorem \ref{theorem 3.27.2}.

\begin{corollary}
                                       \label{corollary 3.26.1}
Assume that there exists a constant $\hat b\in(0,\infty)$
such that, for any $C\in\cC$
\begin{equation}
                                                \label{3.26.2}
\|b\|_{L_{p,q}(C)}\leq\hat b R^{d/p+2/q-1}_{C}.
\end{equation}
Then $\bar b_{\infty}\leq\hat N\hat b$.

\end{corollary}

This follows immediately from
\eqref{3.30.3} with $f=|b|$ and the fact that
a natural modification of \eqref{3.21.80} holds for    
any starting point.

Finally, we need the following.

\begin{lemma}
                                       \label{lemma 3.27.2}
For any $\varepsilon>0$ there exists $\alpha=\alpha(\varepsilon)
>1$ such that for any $R>0$  and $x\in\bR^{d}$ 
$$
\Big(E_{0,x}\Big(\int_{0}^{\tau_{R}}|b(t,x_{t})|\,dt
\Big)^{\alpha}\Big)^{1/\alpha}\leq(1+\varepsilon)\bar b_{\infty}R.
$$
\end{lemma}

Proof. We claim that if $\xi\geq0$, $E\xi\leq A$,
 and $E\xi^{2}\leq 2A^{2}$,
then for any $\varepsilon>0$ there exists $\alpha=\alpha(\varepsilon)
>1$ such that 
$$
E\xi^{\alpha}\leq (1+\varepsilon)A^{\alpha}.
$$
Indeed, by normalizing $\xi$ we may assume that $A=1$. Then
$E\xi^{2}\leq 2$ and for $\alpha\in[1,3/2]$
$$
\frac{d}{d\alpha}E\xi^{\alpha}\leq N,
$$
where $N$ is an absolute constant. This proves the claim.
This also proves the lemma after setting
$$
\xi= \int_{0}^{\tau_{R}}|b(t,x_{t})|\,dt 
$$
and using Remark \ref{remark 3.27.1}. The lemma is proved.

\mysection{ Proof of Theorem \protect\ref{theorem 3.27.2}}
                                           \label{section 4.2.1}

Suppose that \eqref{3.26.20} holds for 
any $R\in(0,\infty)$ and $C\in\cC_{R}$ with $\hat b$
satisfying
$$
\hat N\hat b\leq (2\bar N)^{-1},
$$
where $\bar N$ is taken from Theorem \ref{theorem 8.2.1}
and $\hat N$ is taken from Theorem \ref{theorem 3.21.2}.
For $\lambda\in[0,\infty)$ denote by $x^{\lambda}_{t}$
the diffusion process corresponding to $\lambda b$
in place of $b$ and use the superscript $\lambda$ for other
objects related to $x^{\lambda}_{t}$. Call a $\lambda$
``good'' if  (cf.~Assumption
\ref{assumption 12.18.2}) 
$$
\bar N\bar b^{\lambda}_{\infty}<1,
$$
so that, for $x^{\lambda}_{t}$ in place of $x_{t}$,
 the assertions of Theorem \ref{theorem 3.21.2}
and, hence, \eqref{3.30.3}
 hold true.
 Let $\Lambda$ be the set of good
$\lambda$'s. 
Our claim is that $1\in \Lambda$. Observe that $0\in\Lambda$.

We are going to use the method of continuity proving,
first, that $\Lambda\cap[0,1]$ is closed and, second,
that $\Lambda$ is open to the right
 (and therefore contains points
even beyond $1$). 

If $\lambda_{n}\in \Lambda\cap[0,1]$, $n=1,2,...$,
 converge to $\lambda_{0}$,
then by Corollary \ref{corollary 3.26.1} we have
$\bar b^{\lambda_{n}}_{\infty}\leq \hat N\hat b$, that is
\begin{equation}
                                                \label{3.26.41}
E_{t,x}\int_{0}^{\tau^{\lambda_{n}}_{C}}
 \lambda_{n} |
b(\sft_{s},x^{\lambda_{n}}_{s})|\,ds
\leq \hat N\hat b R
\end{equation}
for any $(t,x)\in \bR^{d+1}$, $R>0$, and $C\in \cC_{R}$,
where $\tau^{\lambda }_{C}$ is the first exit time
of $(\sft_{s},x^{\lambda }_{s})$ from $C$.
By using Girsanov's theorem and Fatou's lemma
one easily shows that \eqref{3.26.41} is also true for $n=0$.
But in that case, $\bar N\bar b^{\lambda_{0}}_{\infty}\leq
\bar N\hat N\hat b \leq 1/2<1$ so that, indeed,
$\Lambda\cap[0,1]$ is closed.

To prove that $\Lambda$ is open to the right,
 first take $\lambda=0$,
$\varepsilon>0$,
$(t,x)\in \bR^{d+1}$, $R>0$, and $C\in \cC_{R}$ and
observe that since $b$ is bounded and has compact support
and $\tau^{ \varepsilon}_{C}\leq R^{2}$, there is a constant $K$ such that
$$
E_{t,x}\int_{0}^{\tau^{ \varepsilon}_{C}}|{  \varepsilon }
b(\sft_{s},x^{ \varepsilon}_{s})|\,ds
\leq \varepsilon KR.
$$
Hence, for $\varepsilon$ small enough
we have $\bar N\bar b^{\varepsilon}_{\infty}<1$,
so that all small $\varepsilon$'s are good.
Next,  take a
$\lambda\in \Lambda \cap(0,1] $,   $\varepsilon>0$,
$(t,x)\in \bR^{d+1}$, $R>0$, and $C\in \cC_{R}$
 and use Girsanov's theorem
to see that
\begin{equation}
                                                \label{3.26.4}
E_{t,x}\int_{0}^{\tau^{\lambda+\varepsilon}_{C}}| (\lambda+\varepsilon) 
b(\sft_{s},x^{\lambda+\varepsilon}_{s})|\,ds
=E_{t,x}e^{\phi(\varepsilon)}\int_{0}^{\tau^{\lambda }_{C}}| (\lambda+\varepsilon) 
b(\sft_{s},x^{\lambda }_{s})|\,ds,
\end{equation}
where  
$$
\phi(\varepsilon)=\varepsilon\int_{0}^{\infty}
 \sigma^{-1}b(\sft_{s},x_{s})\,dw_{s}
-(\varepsilon^{2}/2)
\int_{0}^{\infty}| \sigma^{-1}b(\sft_{s},x_{s})|^{2}\,ds.
$$
 
Recall that $E_{t,x}e^{\phi(\beta\varepsilon)}=1$ for any $\beta$
and observe that for any $\beta>1$
$$
E_{t,x}e^{\beta\phi(\varepsilon)}=
E_{t,x}e^{\phi(\beta\varepsilon)}
\exp\Big((\varepsilon^{2}/2)( \beta^{2} -1)
\int_{0}^{\infty}| \sigma^{-1}b(\sft_{s},x_{s})|^{2}\,ds\Big)
\leq e^{\varepsilon^{2}\beta^{2}K}
$$
since $b$ is bounded and the range of $t$ such that $b(t,\cdot)
\not\equiv 0$ is bounded,
where $K$ is a constant independent of $t,x$.  We use this
and H\"older's inequality to obtain from
\eqref{3.26.4} that
$$
E_{t,x}\int_{0}^{\tau^{\lambda+\varepsilon}_{C}}| (\lambda+\varepsilon) 
b(\sft_{s},x^{\lambda+\varepsilon}_{s})|\,ds
$$
\begin{equation}
                                                 \label{3.27.4}
\leq e^{\varepsilon^{2}\beta KT} 
\Big(E_{t,x}\Big(\int_{0}^{\tau^{\lambda }_{C}}| (\lambda+\varepsilon) 
b(\sft_{s},x^{\lambda }_{s})|\,ds\Big)^{\alpha}\Big)^{1/\alpha},
\end{equation}
where $\alpha=\beta/(\beta-1)$. 
 
Recall that $\lambda$ is good, so that, for any $\varepsilon_{1}
>0$
according to  Lemma \ref{lemma 3.27.2},
for an appropriate choice of $\beta$, the second factor
on the right in \eqref{3.27.4} is less than
$(1+\varepsilon_{1})
 (1+\varepsilon/\lambda) \bar b^{\lambda}_{\infty}R_{C}$,
which by Corollary \ref{corollary 3.26.1} is less than
$(1+\varepsilon_{1}) (1+\varepsilon/\lambda) \hat N\hat b R_{C}$.
Since we can choose
$\varepsilon$ and $\varepsilon_{1}$ arbitrarily, we can make
the left-hand side less than $(3/2)\hat N\hat b R_{C} $.
This shows that $\bar b^{\lambda+\varepsilon}_{\infty}
\leq (3/2)\hat N\hat b  $. Now the condition
 $\hat N\hat b\leq (2\bar N)^{-1}$ implies that $
\bar N \bar b^{\lambda+\varepsilon}_{\infty}<1$, 
so that $\lambda+\varepsilon$ is good
for all small enough $\varepsilon>0$ and this
 brings the proof of the theorem to an end.

\mysection{Proof of Theorem \protect\ref{theorem 3.27.20}}
                                           \label{section 4.2.2}
We take $\hat b=\hat b(d,\delta)$ from 
Theorem \ref{theorem 3.27.2} and 
split the proof into two steps.

{\em Step 1\/}. First we want  to prove that
\eqref{3.30.3} holds if $R\leq R_{0}/2$. To do that
take a smooth $\zeta(x)$ such that
$\zeta=1$ on $B_{R}$, $\zeta=0$ outside $B_{2R}$,
and $0\leq\zeta\leq1$ everywhere. Observe that
 for any $\rho\in(0,\infty)$ and $C\in\cC_{\rho}$,
\begin{equation}
                                                \label{3.29.4}
\|\zeta b\|_{L_{p_{0},q_{0}}(C)}\leq\hat b 
\rho^{d/p_{0}+2/q_{0}-1} .
\end{equation}
Indeed, if $\rho\leq R_{0}$, this follows from the 
assumption of the theorem. However, if $\rho\geq R_{0}$,
then $\rho\geq 2R$ and \eqref{3.29.4} follows from
$$
\|\zeta b\|_{L_{p_{0},q_{0}}(C)}\leq
 \| b\|_{L_{p_{0},q_{0}}(C_{2R})}
\leq
\hat b 
(2R)^{d/p_{0}+2/q_{0}-1}.
$$

After that let $\hat x_{t}$ be the process with drift $\zeta b$.
For this process \eqref{3.30.3} holds for all $R\in(0,\infty)$.
Since the coefficients of $(\sft_{s},x_{s})$ in $\bR\times B_{R}$
coincide with the coefficients
of $(\sft_{s},\hat x_{t})$ and the coefficients are smooth, the distributions
of this processes coincide before they exit from $\bR\times B_{R}$.
Therefore, the left-hand side of \eqref{3.30.3}
does not change if we replace there $x_{t}$ with
$\hat x_{t}$ and we are done with the first step.

{\em Step 2, general $R>R_{0}/2$\/}. 
Applying the same argument, based on the fact that
on the small scale $x_{t}$ behaves like a process with small 
$\bar b_{\infty}$, and using 
Theorem \ref{theorem 1.24.1} we see that with strictly positive
 probability,
  depending only on $d,\delta$, and $R_{0}$,
the process starting at a point $(t_{0},x_{0})$
 will reach
$\Gamma:=[t_{0}+S_{1},t_{0}+S_{2}]\times \bar B_{R_{0}/10}(x_{0}+R_{0}e_{1}/4)$,
where $e_{1}$ is the first basis vector,
$0<S_{1}<S_{2}<\infty$ and the $S_{i}$'s
depend only on $d,\delta$. Repeating this
argument after the process reaches $\Gamma$ 
and taking into account that $R<\infty$, we see that
with probability
$\pi>0$ depending only on $d,\delta$, $R_{0}$, and $R $,
starting from any point in $\bR\times B_{R}$
the process will leave $\bR\times B_{R}$ before time
$T$, where $T$ depends only on $d,\delta$, $R_{0}$, and $R $,
that is
$$
P_{t,x}(\tau'_{R}>T)\leq 1-\pi.
$$
Iterating this inequality we obtain $
P_{t,x}(\tau'_{R}>nT)\leq (1-\pi)^{n}$ for $n=1,2,...$.
This shows, as in the proof of Theorem \ref{theorem 9.1.1},
that to prove the current theorem it suffices to prove that
for   any $R\in(0,\infty)$,  
$x\in\bR^{d}$ 
and Borel $f\geq0$ (notice $\tau_{R}$)
\begin{equation}
                                              \label{3.1.4}
E_{0,x}\int_{0}^{\tau_{R}}
f(t,x_{t})\,dt\leq
\hat N  \|f\|_{L_{p,q} },
\end{equation}
where $\hat N$ depends only on 
$d,\delta ,R_{0}$, and $R$.

Observe that
it suffices to prove
\eqref{3.1.4} only for smooth $f$.  
Fix $\lambda=\lambda(d,\delta)>0$ such that
the right-hand side of \eqref{8.20.1} is less than $1/2$  and introduce
$$
u(t,x)=E_{t,x}\int_{0}^{\tau_{R}}e^{-\lambda s}
f(\sft_{s},x_{s})\,ds.
$$
Then $u$ is a continuous (smooth) nonnegative function 
on $\bar C_{R}$ vanishing on  $\bar C_{R}\setminus (\{0\}\times
B_{R})$
and hence attains   $\max_{\bar C_{R}}u=:M$ at a point
$(t_{0},x_{0})\in C_{R}$. Let $\gamma$ be the first exit time
of $(\sft_{s},x_{s})$ from $C_{R_{0}/2}(t_{0},x_{0})$.
By the strong Markov property
$$
M=u(t_{0},x_{0})=E_{t_{0},x_{0}}e^{-
\lambda(\tau_{R}\wedge\gamma})
u(\sft_{\tau_{R}\wedge\gamma},x_{\tau_{R}\wedge\gamma})
$$
\begin{equation}
                                                \label{3.30.1}
+E_{t_{0},x_{0}}\int_{0}^{\tau_{R}\wedge \gamma}e^{-\lambda s}
f(\sft_{s},x_{s})\,ds.
\end{equation}
Here the second term admits estimating like in
\eqref{3.30.3} by the first step.
The first term is less than
$$
ME_{t_{0},x_{0}}e^{-
\lambda\gamma}I_{\gamma<\tau_{R}}\leq 
ME_{t_{0},x_{0}}e^{-
\lambda\gamma}\leq(1/2)M.
$$
Thus, \eqref{3.30.1} implies that
$$
M\leq N\|f\|_{L_{p,q}}+(1/2)M,\quad M\leq N,
$$
and to finish the proof it only remains to observe that
$$
E_{0,x}\int_{0}^{\tau_{R}}
f(t,x_{t})\,dt\leq e^{\lambda R^{2}}u(0,x).
$$
The theorem is proved.

\mysection{Appendix: a version of Gehring's lemma}
                                     \label{section 3.21.1}

Here we prove the parabolic
 version of the famous Gehring's lemma
stated as Proposition 1.3  in \cite{GS_82}
without proof with the only hint that the proof
is similar to the one given
in the elliptic case in \cite{GM_79}. 
  The author found it quite hard
to make constructions in parabolic case
``similar'' to the elliptic ones given in
\cite{GM_79} and decided to give a complete
proof having a strong probabilistic flavor.
One might think that the only difference
between elliptic and parabolic cases is
different scaling. However, in the elliptic case the doubled cubes  strictly contain 
the original ones and in the parabolic case this is not so.   
Our proof is based on the ideas from \cite{GM_79}   but the organization
of the proof is different. In particular, this allows
us to easily track down the dependence of constants
on $A$ and show that $q$ is a decreasing function of $A$, which was never done before
to the best of the author's knowledge.
If $C=C_{R}(t,x)$ and $\mu>0$ by $\mu C$ we mean
$C_{\mu R}(t,x)$.

\begin{theorem}
                                           \label{theorem 3.9.1}
Let in $C_{R}$ be given a measurable $f(t,x)\geq0$
such that, for some fixed $p,A,B,\mu\in(1,\infty)$ 
satisfying $A\leq B$ and   for all 
 $C\in\cC$  such that $\mu C\subset C_{R}$ we have
$$
\Big(\dashint_{C}f^{p}\,dz\Big)^{1/p}
\leq A\dashint_{\mu C}f\,dz.
$$
Then there exists $q=q(d,p,B)>p$ such that      
$$
\Big(\dashint_{C_{R/4}}f^{q}\,dz\Big)^{1/q}
\leq N(d,p,\mu)A \dashint_{  C_{R/2}}f\,dz.
$$
\end{theorem}

Proof. It is convenient to work with parabolic 
boxes instead of cylinders. For $n= 0,1,...$
and $ k_{0}=0,1,...,2^{2(n+1)  }-1,k_{i}=-2^{n},-2^{n}+1,
...,2^{n}-1$,
 for $i\geq1$, introduce   $D_{k_{0},...,k_{d}}(n)$
 as
$$
[k_{0}2^{-2n},(k_{0}+1)2^{-2n})
\times[k_{1}2^{-n},(k_{1}+1)2^{-n}) 
\times...\times[k_{d}2^{-n},(k_{d}+1)2^{-n}).
$$
We call $2^{-n}$ the size of $D_{k_{0},...,k_{d}}(n)$.
These are dyadic parabolic boxes,
 subsets of $D_{0} :=[0,4)\times[-1,1)^{d}$.
Set $D_{1}=[0,1)\times[-1/2,1/2)^{d}$ and for any
box $D=[S,S+T)\times Q$, where $Q$ is a cube in $\bR^{d}$,
denote   $2D=[S,S+4T)\times 2Q$, where $2Q$
is the concentric cube with twice the side length of $Q$.

Routine arguments show that to prove the theorem,
it suffices to show that there exists 
$q=q(d,p,B)>p$ such that
\begin{equation}
                                       \label{3.9.1}
\Big(\dashint_{D_{1}}f^{q}\,dz\Big)^{1/q}
\leq N(d,p )A \dashint_{ 2D_{0}}f\,dz,
\end{equation}
provided that a nonnegative $f$ is defined in $2D_{0}$ and
\begin{equation}
                                       \label{3.21.5}
\Big(\dashint_{D}f^{p}\,dz\Big)^{1/p}
\leq A\dashint_{2D}f\,dz,
\end{equation}
for any $D=D_{k_{0},...,k_{d}}(n)$ such that
$ D\subset D_{0}$.
 
To proceed in so modified setting,
  for $n\geq 0$ introduce $\Sigma_{n}$ as 
the collection
of $D_{k_{0},...,k_{d}}(n)$. To be consistent with probability
language we add to $\Sigma_{n}$ the empty set.
 Then in the terminology
from \cite{Kr_08}
the family $\{\Sigma_{n}\}$
is a filtration of partitions of $D_{0}$. Observe that
for each $n\geq 0$ and $(t,x)\in D_{0}$ there is only one
element of $\Sigma_{n}$ containing $(t,x)$. We denote
it by $\Gamma_{n}(t,x)$.
Then for each $(t,x)\in D_{0}$ define $\gamma(t,x)$
as the least $n\geq 0$ such that $ 3\Gamma_{n}(t,x)
\subset D_{0}$. Clearly, if $\gamma(t,x)=n$ and $(s,y)
\in \Gamma_{n}(t,x)$, then $\gamma(s,y)=n$. Therefore,
the set $\{(t,x): \gamma(t,x)=n\}$ is the union of
some disjoint elements of $\Sigma_{n}$.
In the terminology
from \cite{Kr_08} this means that $\gamma$ is 
a stopping time relative to the filtration  $\{\Sigma_{n}\}$.

For
each $n\geq 0$ and measurable function $g\geq0$ on $D_{0}$
 one defines the function $g_{|n}$
which on each $\Gamma\in \Sigma_{n}$ equals its average 
over $\Gamma$.

Then for a fixed $\lambda>0$ and $(t,x)\in D_{0}$ we define 
$$
\tau_{\lambda}(t,x)=\inf\{m\geq \gamma(t,x):g_{|m}(t,x)>\lambda\},
\quad(\inf\emptyset:=\infty).
$$
The set $\{\tau_{\lambda}<\infty\}$ is similar to
what one usually
gets by applying the Riesz-Calder\'on-Zygmund decomposition.
However, we are watching the averages of $g$ only on dyadic boxes
where $\gamma$ is constant. Otherwise
we continue in the usual way.

Observe that $D_{0}\cap\{g>\lambda\}\subset  
D_{0}\cap \{\tau_{\lambda}<\infty\}$ (a.e.) because of the Lebesgue
differentiation theorem. 

Next, assume that, for a constant $\bar g$, we have $g_{|\gamma}
\leq \bar g$ and take $\lambda>\bar g$ so that $\tau>\gamma$. Then note that
the set $D_{0}\cap \{ \tau_{\lambda}<\infty\}$
 is either empty
or is the disjoint union
of some nonempty $\Gamma_{i}\in \Sigma_{m_{i}}$, $i=1,2,...$,
on each of which $\tau_{\lambda}=m_{i}$.
Trivially,
$$
\int_{\Gamma_{i}}g\,dz=\int_{\Gamma_{i}}g_{|m_{i}}\,dz
=\int_{\Gamma_{i}}g_{|\tau_{\lambda}}\,dz,
$$
which implies that
$$
\int_{D_{0}}gI_{ \tau_{\lambda}<\infty}\,dz
= \int_{D_{0}}g_{|\tau_{\lambda}}
I_{ \tau_{\lambda}<\infty}\,dz.
$$
Furthermore, on the set $D_{0}\cap \{ \tau_{\lambda}<\infty\}$
we have $g_{|\tau_{\lambda}}>\lambda$, $g_{|\tau_{\lambda}-1}
\leq \lambda$ and, since $g_{|m}\leq 2^{d+2}g_{m-1}$,
we have $g_{|\tau_{\lambda}}\leq\nu^{-1}\lambda$,
where $\nu=2^{-d-2}$. It follows that
$$
\nu\lambda^{-1}\int_{D_{0}}gI_{g>\lambda}\,dz\leq
\nu\lambda^{-1}\int_{D_{0}}gI_{\tau_{\lambda}<\infty}\,dz
=\nu\lambda^{-1}\int_{D_{0}}g_{|\tau_{\lambda}}
I_{\tau_{\lambda}<\infty}\,dz
$$
\begin{equation}
                                               \label{3.21.4}
\leq
|D_{0}\cap\{ \tau_{\lambda}<\infty\}|.
\end{equation}
  We apply this to $g=\phi f^{p}$,
where $\phi(t,x)=[(4-t)^{1/2}\wedge\min_{i}(1-|x^{i}|)]^{d+2}$. As is easy to see
on $D_{0}$ we have
\begin{equation}
                                            \label{3.21.1}
(\phi f^{p})_{|\gamma}\leq N(d)\int_{D_{0}}f^{p}\,dz=:\bar g.
\end{equation}

Next, define $\tilde\Gamma_{1}$ as the largest
(by size)
of the above $\Gamma_{i}$'s and by induction set
$\tilde\Gamma_{i+1}$ to be one of the largest
of $\{\Gamma_{k},k=1,2,...\}\setminus
\{\tilde \Gamma_{k},k=1,2,...,i\}$ such that
its double has no intersection with the doubles
of $\{\tilde \Gamma_{k},k=1,2,...,i\}$. There could be many $\tilde\Gamma_{i }$'s of the same size.  Let
$s_{i}$ denote the size of $\tilde \Gamma_{i}$. We claim that
\begin{equation}
                         \label{4.11.1}
|D_{0}\cap\{\tau_{\lambda}<\infty\}|\leq N(d)
\sum_{i}|\tilde \Gamma_{i}|.
\end{equation}

To prove \eqref{4.11.1} define $\hat\Gamma_{i}$
to be the union of $5\tilde \Gamma_{i}$ and its reflection in its lower base. It turns out that
\begin{equation}
                         \label{4.11.2}
D_{0}\cap\{\tau_{\lambda}<\infty\}\subset
\bigcup_{i}\hat\Gamma_{i}.
\end{equation}
Indeed, if it is not true, then there is a $\Gamma_{i}$,
which is not completely covered by the right-hand side of \eqref{4.11.2}. Let $s$ be the size of $\Gamma_{i}$.
Then there is the largest $k$ such that $s_{k}\geq s $ and $2\Gamma_{i}$ has a nonempty intersection
with at least one of $2\tilde \Gamma_{r}$, $r\geq k$
(because otherwise $\Gamma_{i}\in\{\tilde\Gamma_{r},r\leq k+1\}$). Then, since 
$s_{k}\geq s $, as is easy to see, $\Gamma_{i}
\subset \hat\Gamma_{k}$. This proves \eqref{4.11.2},
which owing to $|\hat\Gamma_{i}|\leq 2\cdot 5^{d+2}
|\tilde\Gamma_{i}|$, implies \eqref{4.11.1}.

Also note that, since $\tau>\gamma$, each of  $\tilde\Gamma^{i}$ 
is a parabolic dyadic box of size $2^{-m_{i}}$
which is the subset of a parabolic dyadic box, say
$\check\Gamma^{j}$,
of size $2^{-k}$, where  $k\leq m_{i}$ is the value
of $\gamma$ on $\check\Gamma^{j}$. It follows by construction
  that
$3\check\Gamma^{j}\subset D_{0}$. In particular, $3\tilde\Gamma^{i}
\subset D_{0}$. Also   the ratio
$\phi(z_{1})/\phi(z_{2})$ is bounded by a constant $N$
as long as $z_{1},z_{2}\in \tilde\Gamma^{i}$.
Therefore,
$$
\lambda|\tilde\Gamma^{i}|^{p}
\leq  
|\tilde\Gamma^{i}|^{p}\dashint_{\tilde \Gamma^{i }}\phi f^{p}\,dz\leq N
|\tilde \Gamma^{i}|^{p}\max_{\tilde \Gamma^{i}}\phi
\dashint_{\tilde \Gamma^{i}}  f^{p}\,dz
$$
$$
\leq
NA^{p} \min_{\tilde \Gamma^{i}}\phi
\Big(\int_{2\tilde\Gamma^{i}}f\,dz\Big)^{p}\leq
NA^{p}  
\Big(\int_{2\tilde\Gamma^{i}}\phi^{1/p}f\,dz\Big)^{p},
$$
$$
|\tilde\Gamma^{i}|\leq N_{1}\frac{A}{ \lambda^{1/p}}
\int_{2\tilde\Gamma^{i}}\phi^{1/p}f\,dz.
$$
One of inconveniences of the last estimate is that we do not have
control of $f$ on $2\tilde\Gamma^{i}$. 
In a similar situation Gehring suggested to sacrifice
some part of what is on the right to be absorbed by the left-hand side
but restrict values of $f$. So following him 
we dominate the right-hand side
by
$$
N_{1}\frac{A}{ \lambda^{1/p}}
\int_{2\tilde\Gamma^{i}}I_{\phi f^{p}>s}\phi^{1/p}f\,dz+
N_{1}\frac{As^{1/p}}{ \lambda^{1/p}}|2\tilde\Gamma^{i}|,
$$
where $s>0$ is arbitrary. For $s= N^{-p}_{2}A ^{-p}\lambda$,
where $N_{2}= N_{1}2^{d+2} $, we get
$$
|\Gamma^{i}|\leq N\frac{ A}{ \lambda^{1/p}}
\int_{2 \tilde\Gamma^{i} }I_{\phi f^{p}>s}\phi^{1/p}f\,dz
$$
and hence, coming back to \eqref{3.21.4}
(and recalling that $2 \tilde\Gamma^{i}$'s are disjoint and $3\tilde\Gamma^{i}
\subset D_{0}$),
for any $\lambda>\bar g$, we obtain
$$
\nu\lambda^{-1}\int_{D_{0}}\phi
f^{p}I_{\phi f^{p}>\lambda}\,dz
\leq  N A\lambda^{-1/p}\int_{  D_{0} }
\phi^{1/p}f I_{\phi f^{p}>N^{-p}_{2}A ^{-p}\lambda}\,dz.
$$
Multiply both sides by $\lambda^{\alpha}$, $\alpha\in(0,1]$, 
and integrate between $\bar g$
and an arbitrary finite $\Lambda>\bar g$ to get
$$
\alpha^{-1}\int_{D_{0}}\phi f^{p} ((\phi f^{p})\wedge \Lambda)^{\alpha}\,dz
-\alpha^{-1}\int_{D_{0}}\phi f^{p} ((\phi f^{p})\wedge \bar g)^{\alpha}\,dz
$$
$$
\leq N (\alpha+1-1/p)^{-1}A \int_{D_{0}}\phi^{1/p}f
  \Big((N_{2}A\phi^{1/p}f)^{p}\wedge \Lambda\Big)^{\alpha+1-1/p}\,dz .
$$

Here
$$
\int_{D_{0}}f^{p} ((\phi f^{p})\wedge \bar g)^{\alpha}\,dz
\leq \bar g^{\alpha}\int_{D_{0}}\phi f^{p}\,dz
\leq N \big(\dashint_{D_{0}} f^{p}\,dz\Big)^{1+\alpha}.
$$
Also 
$$
\phi^{1/p}f \Big((N_{2}A\phi^{1/p}f)^{p}\wedge \Lambda\Big)^{\alpha+1-1/p}
\leq (N_{2}A)^{p(\alpha+1)-1}\phi^{1/p}f((\phi f^{p})\wedge \Lambda)^{\alpha+1-1/p}
$$
$$
\leq (N_{2}A)^{p(\alpha+1)-1}\phi f^{p}((\phi f^{p})\wedge \Lambda)^{\alpha}.
$$
We   conclude that
$$
\int_{D_{0}}\phi f^{p} ((\phi f^{p})\wedge \Lambda)^{\alpha}\,dz
\leq N\Big(\int_{D_{0}}f^{p}\,dz\Big)^{1+\alpha}
$$
$$
+
N_{3}\alpha(\alpha+1-1/p)^{-1}A^{p(\alpha+1)}
\int_{D_{0}}\phi f^{p} ((\phi f^{p})\wedge \Lambda)^{\alpha}\,dz.
$$
Now choose $\alpha\leq1$ so that
$$
N_{3}\alpha(\alpha+1-1/p)^{-1}B^{2p }\leq 1/2.
$$
Then we obtain
$$
\int_{D_{0}}\phi f^{p} ((\phi f^{p})\wedge \Lambda)^{\alpha}\,dz
\leq N\Big(\int_{D_{0}}f^{p}\,dz\Big)^{1+\alpha},
$$
which after sending $\Lambda\to\infty$ and using
\eqref{3.21.5} yields the result with $q=p(1+\alpha)$.
 The theorem is proved.

{\bf Acknowledgment}. The author
thanks the referee for pointing out several
glitches in the paper which resulted in 
correcting and improving the presentation.
Also the author's gratitude is due to
Shuntaro Tsubouchi for pointing out a gap
in the proof of Theorem \ref{theorem 3.9.1},
correcting which led to a shorter proof.


\begin{thebibliography}{mm}

\bibitem{Al_60} A. D. Aleksandrov, {\em
Certain estimates for the Dirichlet problem\/},
Dokl. Akad. Nauk SSSR, Vol. 134 (1960),
 1001--1004 (Russian); translated as
Soviet Math. Dokl., Vol. 1 (1961), 1151--1154.

\bibitem{Al_63} A. D. Aleksandrov,
{\em Uniqueness conditions and estimates for
the solution of the Dirichlet
problem\/}, Vestnik Leningrad. Univ., Vol. 18 (1963), No. 3,
5-29 in Russian;
English translation in
Amer. Mat. Soc. Transl., Vol. 68 (1968), No. 2, 89-119.

\bibitem{BFGM_19} L. Beck, F. Flandoli, M. Gubinelli, and M. Maurelli,
{\em Stochastic ODEs and stochastic linear PDEs with
critical drift: regularity, duality and uniqueness\/},
Electron. J. Probab., Vol. 24 (2019), No. 136, 1--72.



\bibitem{Ca_95} X. Cabr\'e, 
{\em On the Alexandroff-Bakelman-Pucci estimate and the reversed
H\"older inequality for solutions of elliptic 
and parabolic equations\/}, Comm. Pure
Appl. Math., 48 (1995), 539--570.


\bibitem{CKS_00} M. G. Crandall, M. Kocan, and A. \'Swi{\c e}ch, {\em
$L^p$-theory for fully nonlinear uniformly parabolic equations\/},
 Comm. Partial Differential Equations, Vol. 25  (2000),
 No. 11-12, 1997--2053.

\bibitem{DK_21}  Hongjie Dong and N.V. Krylov,
{\em
Aleksandrov's estimates for elliptic equations with drift in a Morrey
  spaces containing $L_{d}$\/},  http://arxiv.org/abs/2103.03955

\bibitem{FS_84} E.B. Fabes and D.W. Stroock,
{\em The $L^{p}$-integrability of Green's functions
and fundamental solutions for elliptic
and parabolic equations\/}, Duke Math. J., Vol. 51
(1984), No. 4, 997--1016.

\bibitem{Fo_98} K. Fok, {\em A nonlinear Fabes-Stroock result\/},
Comm. Partial Differential Equations, 23 (1998), No. 5-6, 967--983.

\bibitem{GM_79} M. Giaquinta and G. Modica,
{\em Regularity results for some classes of higher
order non linear elliptic systems\/}, J. Reine Angew. Math.,
Vol. 311(312) (1979), 145-169.

\bibitem{GS_82} M.Giaquinta and M. Struwe,
{\em On the partial regularity of weak solutions 
of nonlinear parabolic systems\/}, 
Mathematische Zeitschrift, Vol. 179 (1982), 437-451. 

\bibitem{Kr_74} N.V. Krylov, 
{\em  Some estimates for the density of  distribution
  of a
 stochastic
 integral\/}, Izvestiya Akademii Nauk SSSR, seriya matematicheskaya,
Vol. 38 (1974), No. 1,  228--248 in Russian; English translation
in Math. USSR
Izvestija,   Vol. 8  (1974),  No. 1,   233--254.

\bibitem{Kr_77}  N.V. Krylov,  ``Controlled diffusion processes'',
Nauka, Moscow,  1977 in Russian; English  transl. 
   Springer,
1980.

\bibitem{Kr_86} N.V. Krylov,
  {\em  On  estimates of the maximum of a solution of
a parabolic equation and  estimates of the distribution of a
semimartingale}, Matematicheski Sbornik, Vol. 130, No. 2  (1986),  207--221
in Russian, English translation is
 Math. USSR Sbornik,  Vol. 58 (1987), No. 1,  207--222.

\bibitem{Kr_08} N.V. Krylov,
``Lectures on elliptic and parabolic equations
in Sobolev spaces", Amer.
Math. Soc., Providence, RI, 2008.


\bibitem{Kr_19_1} N.V. Krylov, {\em
On stochastic equations with drift in $L_{d}$\/},
 Annals of Prob.,   Vol. 49 (2021), No. 5, 2371--2398.

\bibitem{Kr_20} N.V. Krylov,  {\em
Linear and fully nonlinear elliptic equations with $L_{d}$-drift\/},
  Comm. PDE, Vol. 45 (2020), No.  12, 1778--1798.

\bibitem{Kr_20_2} N.V. Krylov, {\em
On time inhomogeneous stochastic It\^o equations 
with drift in $L_{d+1}$\/},  
Ukrains'kyi Matematychnyi Zhurnal, Vol. 72 (2020), No. 9,  1232--1253.

\bibitem{Kr_21_1}  N.V. Krylov, {\em
On potentials of It\^o's processes with drift in $L_{d+1}$\/},\\
  http://arxiv.org/abs/2102.10694

\bibitem{Kr_21_2}  N.V. Krylov, {\em
On diffusion processes with drift in $L_{d+1}$\/},\\
  http://arxiv.org/abs/2102.11465

\bibitem{Kr_21_3} N.V. Krylov, {\em
On the heat equation  
with drift in $L_{d+1}$\/},\\
arXiv:2101.00119

  \bibitem{Kr_21b}  N.V. Krylov,  {\em
Elliptic equations with VMO a,
  b$\,\in L_{d}$, and c$\,\in L_{d/2}$\/}, Trans. Amer. Math.
Sci., Vol. 374  (2021), No. 4, 2805-2822. 

 

\bibitem{Na_15} A.I. Nazarov, 
{\em
Interpolation of linear spaces and estimates for 
the maximum of a solution for parabolic equations\/}, Partial differential
equations, Akad.  Nauk SSSR Sibirsk. Otdel., Inst. Mat., No\-vo\-si\-birsk,
1987, 50--72 in Russian; translated into English 
as  {\em
On the maximum principle for parabolic equations
with unbounded coefficients\/},
https:// arxiv.org/abs/1507.05232

\bibitem{NU_85} A.I. Nazarov and N.N. Ural'tseva, {\em
Convex-monotone hulls and an estimate of the maximum of the solution of
a parabolic equation\/}, Boundary value problems of
mathematical physics and related problems 
in the theory of functions,
No. 17,  Zap. Nauchn. Sem. Leningrad. Otdel. Mat. Inst. Steklov. 
(LOMI) Vol. 147 
(1985), 95--109, in Russian, English translation in
 Journal of Soviet
Mathematics Vol. 37 (1987), 851--859

\bibitem{RZ_21_1} M. Roeckner and Guohuan Zhao,
{\em SDEs with critical time dependent drifts: weak solutions\/},
 	arXiv:2012.04161

\bibitem{RZ_21_2} M. Roeckner and Guohuan Zhao,
{\em SDEs with critical time dependent drifts: strong solutions\/},
arXiv:2103.05803

\end{thebibliography}
\end{document}